\documentclass[]{article}

\usepackage{amsmath}
\usepackage{amssymb}
\usepackage[ansinew]{inputenc}

\def\Nset{\mathbb{N}}

\def\Rset{\mathbb{R}}

\newtheorem{theorem}{Theorem}
\newtheorem{lemma}{Lemma}
\newtheorem{example}{Example}
\newtheorem{cor}{Corollary}
\newtheorem{defi}{Definition}

\newtheorem{prop}{Proposition}
\newcommand{\proof}{\bfseries Proof: $\quad$\mdseries}
\newcommand{\finishproof}{\begin{flushright} $\Box$ \end{flushright}}

\date{\today}

\begin{document}

\title{Asymptotic optimality of isoperimetric constants with respect to $L^{2}(\pi)$-spectral gaps}
\author{Achim Wübker\thanks{e-mail: awuebker@mathematik.uni-osnabrueck.de}\\Institute of Mathematics, University of Osnabrück\\Albrechtstraße. 28 a, 49076 Osnabrück}

\maketitle

\begin{abstract}
In this paper we investigate the existence of $L^{2}(\pi)$-spectral gaps for $\pi$-irreducible, positive recurrent Markov chains on general state space. We obtain necessary and sufficient conditions for 
the existence of $L^{2}(\pi)$-spectral gaps in terms of a sequence of isoperimetric constants and
establish their asymptotic behavior. It turns out that in some cases the spectral gap can be understood in terms of convergence of an induced probability flow to the uniform flow. The obtained theorems can be interpreted as mixing results and yield sharp estimates for the spectral gap
of some Markov chains.

\end{abstract}

\section{Introduction}
Throughout this paper we consider a time homogeneous and time discrete positive recurrent Markov chain $\xi_1,\xi_2,\ldots$
on an arbitrary state space $(\Omega,\mathcal{F})$ with transition kernel $p(\cdot,\cdot)$ and 
uniquely determined invariant measure $\pi$. It is assumed that the $\sigma$-Algebra $\mathcal{F}$ is countably generated. 
The main result of this paper is the precise characterization of the size of the spectral gap in terms of a sequence of isoperimetric 
constants for normal Markov operators. This characterization has many consequences, for example it will be applied  
to obtain sharp upper bounds for the $L^2(\pi)$-spectral gap for the $M/M/1$-queue, the hypercube and the star. Moreover, the theorem is of 
theoretical value. On the one hand it is used to generalize a result due to Roberts and Tweedie \cite{roberts1} concerning the connection between
geometric ergodicity and the spectral gap property (SGP). On the other hand, it is applied to obtain a functional analytic
inequality between the Banach space of complex valued measurable functions $\mathcal{B}(\Omega,\mathcal{F},||\cdot||_{\infty})$ and the Hilbert
space $L^2$.  \\
The paper is organized as follows: First we introduce the main notions and notations needed throughout this paper, especially a certain sequence of isoperimetric constants, which is used to provide a classification of those Markov chains having the (SGP) (section \ref{verygeneral}).    
In order to obtain more precise information about the size of the spectral gap, we shrink the class of considered Markov chains to 
those corresponding essentially to normal Markov operators on $L^2$. Here we state and prove our main results (section \ref{sectioniso3}). These results are applied in section \ref{apply}.
Finally, it is shown how the results can be ranged in the theory of mixing for Markov chains (section \ref{mixing1}). The proofs, that are not 
presented so far, can be found in section \ref{proofs}.

We say that $P$ has an $L^2(\pi)$-spectral gap if   

\[\lim_{n\rightarrow\infty}\sup_{f\in L_{0,1}^{2}(\pi)}||P^{n}f||_{2}^{\frac{1}{n}}=:\rho<1,\]
where 
\[
Pf(x):=\int_{\Omega}f(y)P(x,dy),\,\,f\in L^{2}(\pi)
\]
and $L_{0,1}^{2}(\pi):=\{f\in L^{2}(\pi):\int_{\Omega}f(x)\pi(dx)=0,\,\int_{\Omega}f(x)^{2}\pi(dx)=1\}$.
The spectral gap is defined to be $1-\rho$.
Let $L_{0}^{2}(\pi):=\{f\in L^{2}(\pi):\int_{\Omega}f(x)\pi(dx)=0\}$. If not stated otherwise, $P$ is regarded as
an operator on the Hilbert-space $L_{0}^{2}(\pi)$. For reversible Markov chains, i.e 
\[\pi(dx)p(x,dy)=\pi(dy)p(y,dx),\]
there is another useful description for $\rho$, namely
\begin{equation}\label{Rayleigh-Ritz_principle}
||P||_{L_{0}^{2}(\pi)}=\sup_{\lambda\in\sigma(P)}\{|\lambda|:\lambda\not=1\}=\rho=\sup_{f\in L_{0,1}(\pi)}|\langle{f,Pf}\rangle_{\pi}|.
\end{equation}

The right hand side of equation (\ref{Rayleigh-Ritz_principle}) is known as the Rayleigh-Ritz principle \cite{lawler}.
Note that for reversible Markov chains 
the associated Markov operator $P$ is selfadjoint and  that the
spectrum $\sigma(P)$ is always contained in $[-1,1]$.\\
\begin{defi}\label{geoergodic}
A Markov chain $\xi_1,\xi_2,\ldots$ is called $\pi$-a.s. geometrically ergodic if there exists $\delta<1$, $C_{x}>0$ such that 
\begin{equation}\label{geometrischergodisch}
||p^{n}(x,\cdot)-\pi(\cdot)||_{TV}\le C_{x}\delta^{n}\,\,\,\,\pi\,-a.s.,
\end{equation}
\end{defi}
where $||\cdot||_{TV}$ denotes the total-variation-norm. In the following we will assume that $\delta$ is chosen optimal, i.e. $\delta$ is as small as
possible satisfying (\ref{geometrischergodisch}). This $\delta$ will be called the optimal rate of geometric convergence (ORGC).
We have the following alternative characterization of geometric ergodicity,
which can be found in
\cite{nummelin2} or \cite{nummelin}, Theorem 6.14 (iii). 
\begin{theorem}\label{habbatu}
Let $\xi_1,\xi_2,\ldots$ be an ergodic Markov chain on a probability space with countable generated $\sigma$-algebra $\mathcal{F}$. Then we have the following equivalence:
\begin{enumerate}
\item\label{jajaja}
\begin{equation}
\exists \delta_{x}<1, C_{x}>0: \,\,\,\,||p^{n}(x,\cdot)-\pi(\cdot)||_{V}\le C_{x}\delta^{n}\,\pi-a.s.
\end{equation}
\item\label{muffi0}
\begin{equation}\exists \delta<1, C>0, \mbox{ such that }||\,\,||p^{n}(\cdot,\cdot)-\pi||_{V}||_{1}\le C \delta^{n}.
\end{equation}
\end{enumerate}
\end{theorem}

There is a close relation between the (SGP) and convergence in the total variation-norm. 
Doeblin \cite{doeblin2} and Doob \cite{doob} had shown that the convergence in (\ref{geometrischergodisch}) holds uniformly 
in $x$ (i.e. $C(x)=C$ $\forall x\in\Omega$) if and only if  
\[\lim_{n\rightarrow\infty}\sup_{f\in \mathcal{B}_{\infty}}||(P^{n}-P_1)f||_{\infty}^{\frac{1}{n}}<1,\]
where $\mathcal{B}_{\infty}$ denotes the space of all bounded, $\mathcal{F}$-measurable complex valued functions and $P_1(f)=E_{\pi}f=\int_{\Omega}f(y)\pi(dy)$.
Using the Riesz-Thorin-interpolation theorem one can show that this implies the existence of an $L^{2}(\pi)$-spectral gap. Actually, we will obtain this result with some additional refinements as a consequence of the theory developed in this paper.

The following families of isoperimetric constants will play a crucial role during this work: 
\[\label{isoneu}
k_{n}:=\inf_{A\in\mathcal{F}}k_{n}(A),\quad k_{n}(A):=\frac{1}{\pi(A)\pi(A^{c})}\int_{A}p^{n}(x,A^{c})\pi(dx), \,\,n\in\Nset,
\]
and
\[k_{P^{\ast^{n}}P^{n}}:=\inf_{A\in\mathcal{F}}k_{P^{\ast^{n}}P^{n}}(A):=\inf_{A\in\mathcal{F}}\frac{1}{\pi(A)\pi(A^{c})}\int_{A}P^{\ast^ {n}}P^{n}1_{A^{c}}(x)\pi(dx),\]
where $P^{\ast}$ is the adjoint operator of $P$ considered on $L^{2}(\pi)$.\\
In \cite{wuebker} it was shown that the following sequence of constants can be seen as an appropriate way to measure the periodic behavior of a Markov chain $\xi_1,\xi_{n+1},\ldots$. 
\[K_n:=\sup_{A\in\mathcal{F}}k_{n}(A), \,\,n\in\Nset.\]

\section{Bounds for the spectral gap in terms of isoperimetric constants in the general case}\label{verygeneral}
Let us present an upper bound for the spectral radius of $P$ in terms of the isoperimetric quantities $k_n$ and $k_{P^{\ast^n}P^n}$. 
\begin{prop}\label{noetig}
Let $\xi_1, \xi_2, \ldots$ be a stationary Markov chain with invariant measure $\pi$, transition kernel
$p(\cdot,\cdot)$ and state space $(\Omega,\mathcal{F})$. Then the condition
\begin{equation}\label{wichtigevoraussetzung}
0<\inf_{n\in\Nset}k_{n}=\inf_{n\in\Nset}\inf_{A\in\mathcal{F}}k_{n}(A).
\end{equation}
is necessary for the (SGP). More precisely, we obtain 
the following estimate:
\begin{equation}\label{wichtigevoraussetzung2}
||P^{n}||^{2}_{2}\ge\sup_{A\in\mathcal{F}}\frac{1}{2}\left(\left(1-k_n(A)\right)^{2}+2\pi(A)\pi(A^{c})\left(1-k_{P^{\ast^n}P^n}(A)\right)\right).
\end{equation}
\end{prop}
The necessity of (\ref{wichtigevoraussetzung}) follows immediately from (\ref{wichtigevoraussetzung2}), since if (\ref{wichtigevoraussetzung}) is not
satisfied, the right hand side in (\ref{wichtigevoraussetzung2}) is bounded away from zero for all $n\in\mathbb{N}$, so $||P^{n}||^{2}_{2}$ is, which
imply that the (SGP) cannot hold. The other direction of the proof is postponed to section \ref{proofs}.

It is natural to ask whether the condition (\ref{wichtigevoraussetzung}) is also
sufficient. In the following, we will show that this is actually the case. In fact, we will see that an apparently weaker condition than (\ref{wichtigevoraussetzung}) turns out to be equivalent to the (SGP). To show this, we will need the following result due to 
Lawler and Sokal \cite{lawler}, which says that

\begin{equation}\label{soiki}
\Re(\sigma(Id-P))\ge \frac{\kappa}{8}k^{2},
\end{equation}
with 
\begin{equation}\label{kappageben}
\kappa=\inf_{\mathcal{D}}\sup_{c\in\Rset}\frac{E\left(|(X+c)^{2}-(Y+c)^{2}|\right)}{E((X+c)^{2})}
\end{equation}
and $\mathcal{D}$ denoting the set of all possible distributions of i.i.d random variables $(X,Y)$ with variance 1 and where $\Re(\sigma(Id-P))$ means the real part of the spectrum of $(Id-P)$.

Additionally, Lawler and Sokal \cite{lawler} had shown that $\kappa\ge 1$.

If we apply the spectral mapping theorem (see e.g. \cite{werner}) to (\ref{soiki}) with the function $f=1-z$, we obtain
\[\Re(\sigma (P))\le 1-\frac{\kappa}{8}k^2.\]
Moreover,
\[\Re(\sigma(P^{n}))\le 1-\frac{\kappa}{8}k_{n}^2.\]
Now let us state our first result:

\begin{theorem}\label{spaet}
Let $\xi_1, \xi_2,\ldots$ be a Markov chain with state space $(\Omega,\mathcal{F},\pi)$ and associated Markov 
operator $P$. Then the following statements are equivalent:
\begin{enumerate}
\item
$P$ has an $L^{2}(\pi)$-spectral gap.
\item
\begin{equation}\label{zwo}
\mathcal{M}:=\left\{\epsilon>0:\,\,\,k_n>\epsilon\,\,\,\,\,\forall n\le\left[\frac{2\pi}{\arccos(1-\frac{\kappa}{16}\epsilon^{2})}\right]+1\right\}\not=\emptyset.
\end{equation}
\end{enumerate}
If condition (\ref{zwo}) holds true, we obtain
\begin{equation}\label{dro}
\sigma(P-P_1)\subset B_{r_0}(0),\,\,r_{0}:=\inf_{\epsilon\in\mathcal{M}}(1-\frac{\kappa}{16}\epsilon^{2})^{1/([\frac{2\pi}{\arccos(1-\frac{\kappa}{16}\epsilon^{2})}]+1)},
\end{equation}
where $B_{r}(0)$ is the ball of radius r with center $0$. 
\end{theorem}
The necessity of (\ref{zwo}) follows immediately from Proposition \ref{noetig}, so the interesting part of the theorem is the sufficiency. 
The theorem says that it is enough to check only finitely many $n$ with $k_n\ge\epsilon$ for establishing the (SGP), where this number again depends on $\epsilon$.\\
Let us study the asymptotic behavior of (\ref{zwo}) and (\ref{dro}) for $\epsilon\rightarrow 0$. First, note that
\begin{equation}\label{appearing}
\frac{2\pi}{\arccos(1-\frac{\kappa}{16}\epsilon^{2})}\sim\frac{8\pi}{\sqrt{2\kappa}}\frac{1}{\epsilon}
\end{equation} 
(in the sense $a_{\epsilon}\sim b_{\epsilon}\Longleftrightarrow \lim_{\epsilon\rightarrow 0}\frac{a_\epsilon}{b_\epsilon}$=1).
This means that the number of $n$ with $k_n>\epsilon$ we have to check in order to apply Theorem \ref{spaet} depends anti-proportional on $\epsilon$.\\

Now we will show that if we do not add further assumptions, this is in a certain sense the weakest sufficient condition (and as shown necessary) for the existence of an $L^{2}(\pi)$ spectral gap in terms of isoperimetric constants $k_n$.
To see this, let $p$ be a prime number and consider the deterministic walk on the circle $\mathbb{Z}/p\mathbb{Z}$, i.e. 
\[\mathbb{Z}/p\mathbb{Z}\ni i\rightarrow i+1\in\mathbb{Z}/p\mathbb{Z}.\]
It is readily seen that this Markov chain is positive recurrent. The invariant distribution can be easily established to
be the uniform distribution. Furthermore, one can check that 
\[k_i>\frac{1}{p}\,\,\,\forall\,\,i<\frac{1}{1/p},\]
but 
\[k_{p}=0.\]
This holds true for all prime numbers $p\in\mathbb{N}$.
So if we set $\epsilon=\frac{1}{p}$, we see that $constant *\frac{1}{\epsilon}$ as suggested in Theorem \ref{spaet} is the right order of magnitude of $k_i$'s to check and this cannot be improved. What certainly can be improved is the 
constant $\frac{8\pi}{\sqrt{2\kappa}}$ appearing in (\ref{appearing}) and it would be interesting to do this, because it would yield sharper bounds for the spectral gap of the Markov chain.\\

Let us refine the asymptotic results obtained in (\ref{appearing}). 
Using Taylor expansion
the right hand side in (\ref{dro}) can be written as
\[\exp(-\frac{\sqrt{2\kappa^3}}{128\pi}\epsilon^3 +O(\epsilon^5)).\]
From these calculations, we obtain
\begin{cor}\label{kubisch}
Assume that the Markov chain $\xi_1,\xi_2,\ldots$ has the (SGP). Then the spectral gap $1-\rho$ can be estimated by
\begin{equation}
1-\rho\ge \frac{\sqrt{2\kappa^3}}{128\pi}k_{-}^3 -O(k_{-}^5)),
\end{equation}
where $k_{-}:=\min\{k_{n},n\in\mathbb{N}\}$.
\end{cor}
Corollary \ref{kubisch} may be interesting for Markov chains with very small spectral gaps. If it is
possible to determine $n_0$ with $k_{n_0}=k_{-}$, especially if $k_1=k_{-}$, Corollary \ref{kubisch} can be used to 
estimate the spectral gap $1-\rho$ from below. For example, if $P$ is a positive and self-adjoint operator, using the spectral representation theorem for $P$ and the technique presented in the next section one can show that $k_1=k_{-}$. Unfortunately, this case is not very interesting, because Lawler and Sokal's result can be directly applied.
Intuitively, monotonie of the $k_n$ is close related
to the property that the Markov chain has at most a "`weak cyclic behavior"'.
In general it seems to be difficult to present
conditions which ensure the monotonicity of the $k_n$. Even for reversible and not positive chains we can give only 
a partial satisfactorily answer, namely, in this case it can be established that $k_{2n}$ is monotonic increasing in $n$. The general problems arise from the difficulty to obtain precise information about the spectral measures
associated to the family of indicator functions. 

Now let us show how Theorem \ref{spaet} can be applied to prove sufficient conditions for the spectral gap property. 
\begin{cor}\label{hidoe}
Assume that the Markov chain $\xi_1,\xi_2,\ldots$ satisfies (\ref{geometrischergodisch}) with \\
$C(x)=C$ independent of $x$. Then $P$ has an $L^{2}(\pi)$-spectral gap.
\end{cor}
\proof

\begin{eqnarray}
k_{n}(A)&=&\frac{1}{\pi(A)\pi(A^{c})}\int_{A}(p^{n}(x,A^{c})-\pi(A^{c}))\pi(dx)+1\nonumber\\
&\ge&1-\frac{2}{\pi(A)}\int_{A}|p^{n}(x,A^{c})-\pi(A^{c})|\pi(dx)\nonumber\\
&\ge& 1-2C\delta^{n}\ge \frac{1}{2},\,\,n\ge n_{0},\,\,n_{0}\mbox{ sufficiently large},
\end{eqnarray}
where $\delta$ is given in (\ref{geometrischergodisch}).

In order to complete the proof, we need
\begin{lemma}\label{dasbraucheichnoch}
Assume that there exists $n_{0}\in\Nset,\epsilon>0$, such that $k_n\ge\epsilon\,\,\,\forall n\ge n_{0}$. Then we have:
\[k_i\ge\frac{\epsilon}{n_0}\,\forall i\in\Nset.\]
\end{lemma}
With Lemma \ref{dasbraucheichnoch} it follows that
\[k_n\ge\frac{1}{2n_0}\,\,\forall n\in\Nset.\]
The claim follows by applying Theorem \ref{spaet}.
\finishproof

Since the assumption in Corollary \ref{hidoe} is equivalent to the classical Doeblin condition (see \cite{tweedie}),
we have just shown the well-known fact that the classical Doeblin condition implies the existence of an $L^{2}(\pi)$-spectral gap.
Later we will improve this result by comparing rates of convergence on different Banach-spaces.

Moreover, from Theorem \ref{spaet} we easily derive the following sufficient condition for the existence of an $L^{2}(\pi)$-spectral gap:
\begin{cor}\label{witiko}
Assume that condition (\ref{zwo}) is fulfilled. Moreover, assume that there exists $\epsilon>0$ such that
\begin{equation}
\limsup_{n\rightarrow\infty}\sup_{A\in\mathcal{F}:\pi(A)\le\frac{1}{2}}\frac{1}{\pi(A^{c})}\int_{A}|\frac{p^{n}(x,A)}{\pi(A)}-1|\pi(dx)\le 1-\epsilon. 
\end{equation}
Then $P$ has an $L^{2}(\pi)$-spectral gap.
\end{cor}

\proof
It is easy to show that
\begin{eqnarray}
k_{n}(A)&\ge& 1-\frac{1}{\pi(A)\pi(A^{c})}\int_{A}|p^{n}(x,A^{c})-\pi(A^{c})|\pi(dx)\nonumber\\
&=&1-\frac{1}{\pi(A^{c})}\int_{A}|\frac{p^{n}(x,A)}{\pi(A)}-1|\pi(dx)\nonumber\\
&\ge&\frac{\epsilon}{2}.
\end{eqnarray}
Now the claim follows from Theorem \ref{spaet} and Lemma \ref{dasbraucheichnoch}.
\finishproof

Until now, our approach provides only a few information about the size of the spectral gap. In order to obtain more precise
estimates for the spectral gap,

we should analyze the asymptotic behavior of $1-k_n$ and $1-k_{P^{\ast^n}P^n}$.

\section{Convergence rates of $k_{n}$ and $k_{P^{\ast^{n}}P^{n}}$}\label{sectioniso3}
In this section we will present the main theorem of this paper, namely the fact that for normal Markov operators we get a precise formulation 
for the spectral gap of $P$ in terms of $1-k_{P^{\ast^n}P^n}$. This result is surprising and has important consequences as will be shown later on.\\
In order to deduce the interesting properties of $k_{n}$ and $k_{P^{\ast^{n}}P^{n}}$, we should mention the 
following elementary but important

\begin{lemma}\label{monotonie}
The sequence $k_{P^{\ast^{n}}P^{n}}(A)$ is monotonic increasing in $n$  and we have 
\begin{equation}\label{equation_monotonie1}
k_{P^{\ast^{n}}P^{n}}(A)\le 1\,\,\, \forall n\in\Nset, \,\,A\in\mathcal{F}.
\end{equation}
\end{lemma}\label{lawmon}

Let us return to Proposition \ref{noetig}. From (\ref{wichtigevoraussetzung2}) it can be immediately obtained
that 
\begin{equation}\label{onlyestimate}
\rho:=\lim_{n\rightarrow\infty}||P^{n}||^{\frac{1}{n}}_2\ge \lim_{n\rightarrow\infty}\max(1-k_n,K_n -1)^{\frac{1}{n}}.
\end{equation}
So the (SGP) implies the geometric convergence of $K_n$ and $k_n$ to $1$, which is equivalent to the geometric convergence of $K_n-k_n$ to zero. 
Indeed, it turns out that this rate of convergence to zero is closely related to the spectral radius $\rho(P)$. To see this,  
first note that (\ref{onlyestimate}) is based on (\ref{wichtigevoraussetzung2}), which sometimes turns out to be not as good as possible from an asymptotic point of view. The following equality, which can be derived by similar calculations as in the proof of 
Proposition \ref{noetig}, indicates that for investigating the spectral properties of $P$, one should have a closer look at the isoperimetric constants $k_{P^{\ast^n}P^n}$:
\begin{equation}\label{evenequality1}
||P^{n}f_{0}||_{2}=\sqrt{1-k_{P^{\ast^n}P^n}},
\end{equation}
with $f_{0}$ as in (\ref{fnot})
so that 
\begin{equation}\label{evenequality2}
\rho\ge \lim_{n\rightarrow\infty}\left(\sqrt{1-k_{P^{\ast^n}P^n}}\right)^{\frac{1}{n}}
\end{equation}
Note that this inequality is a lower bound for the spectral radius of $P$ and hence yield to an upper bound of the spectral gap. We will see that for normal operator actually equality hold in (\ref{evenequality2}).
Before we are going to prove this, let us compare $k_n$ to $k_{P^{\ast^n}P^n}$.\\
Since we have equality in (\ref{evenequality1}) and the estimate in (\ref{onlyestimate}) is also derived by considering $||P^{n}f_{0}||_{2}$, the right hand side in (\ref{evenequality2}) must be at least as large as the right hand side in (\ref{onlyestimate}). 
In general, we have the following relationship between $k_n$ and $k_{P^{\ast^n}P^n}$: 

\begin{lemma}\label{bani}
Let $\xi_1,\xi_2,\ldots$ be a Markov chain with invariant measure $\pi$. Then we have the following inequality
\begin{equation}
1-k_{n}\le\sqrt{2}\sqrt{1-k_{P^{\ast^n}P^n}}.
\end{equation}
\end{lemma}

Having established an upper bound, it would be nice to achieve a lower bound for the spectral gap in terms of 
isoperimetric constants. For this reason let us consider an arbitrary indicator function $g$ in $L^{2}_{0,1}(\pi)$, i.e.
\[g=\sum_{i=1}^{l}\alpha_i\,1_{A_i},\]
where $A_i\in\mathcal{F}$ such that $\bigcup_{i=1}^{l_{g}}A_i=\Omega$, $A_i\cap A_j=\emptyset$ for $i\not=j$, $i,j\in\{1,,\ldots,l_{g}\}$.
Then we have
\begin{eqnarray}\label{laterneed}
||P^{n}g||^{2}&=&\int_{\Omega}\left(\sum_{i=1}^{l_g}\alpha_i p^{n}(x,A_i)\right)^{2}\pi(dx)\nonumber\\
&\stackrel{Jensen}{\le}&\int_{\Omega}\sum_{i=1}^{l}\frac{1}{\pi(A_i)}\int_{\Omega}\left(p^{n}(x,A_i)-\pi(A_i)\right)^{2}\pi(dx)\nonumber\\
&\stackrel{(\ref{adjungrate})}{=}& \sum_{i=1}^{l}\pi(A_{i}^{c})(1-k_{P^{\ast^n}P^n}(A_i))\le l_{g}(1-k_{P^{\ast^{n}}P^{n}}).
\end{eqnarray}
So for all indicator functions $g$ we obtain
\begin{equation}\label{notpas}
\lim_{n\rightarrow\infty}||P^{n}g||^{\frac{1}{n}}\le\lim_{n\rightarrow\infty}\left(\sqrt{1-k_{P^{\ast^n}P^n}}\right)^{\frac{1}{n}}.
\end{equation}
Note that the convergence in (\ref{notpas}) is not uniform with respect to the class of step functions, so we cannot pass to 
limit in order to obtain this estimation for all $f\in L^{2}_{0}(\pi)$. There are mainly two reasonable ways
to continue. First, one could consider only such functions $f\in L^{2}_{0}(\pi)$ which can be sufficiently fast 
approximate by step-functions in a certain sense. For the second way, which we will follow here, one can put some
additional assumptions on the operator $P$. For this reason let us define

\begin{defi}
We call a sequence of operators $(A_n)_{n\in\mathbb{N}}$ positive and quasi selfadjoint generated, if there exists a positive and selfadjoint operator $Q$, constants $C>0$, $q<1$ such that for all indicator functions $g\in L^{2}_{0,1}(\pi)$ 
\begin{equation}\label{q}
||(A_{n}-Q^n) g||\le C q^n\,\,\,\,\,\forall n\in\mathbb{N}. 
\end{equation} 
\end{defi} 

Now let us state the main theorem of this paper:


\begin{theorem}\label{assumdiscuss}
Assume that the sequence $(P^{\ast^n}P^n)_{n\in\mathbb{N}}$ associated to the Markov chain $\xi_1,\xi_2,\ldots$ is quasi
selfadjoint. Then we have for all $f\in L^{2}_{0}(\pi)$
\begin{equation}\label{pas}
\lim_{n\rightarrow\infty}||P^{n}f||^{\frac{1}{n}}_{2}\le\max\left(\lim_{n\rightarrow\infty}\left(\sqrt{1-k_{P^{\ast^n}P^n}}\right)^{\frac{1}{n}},q\right),
\end{equation}
where $q$ is as in (\ref{q}).
\end{theorem}

\proof
For the proof we need the spectral theorem for selfadjoint operators (see e.g. \cite{werner}).
For all $g\in L^{2}_{0,1}(\pi)$ we have 
\begin{equation}\label{eins}
||P^{n}g||^{2}=\langle{Q^{n}g,g}\rangle_{\pi}+\langle{(P^{\ast^n}P^n -Q^n)g,g}\rangle_{\pi}
\end{equation}
Since $Q$ is selfadjoint, we obtain
\begin{eqnarray}\label{zwei}
\langle{Q^{n}g,g}\rangle_{\pi}&=&\int_{\sigma(Q)}\lambda^{n}\langle{d\,E_{\lambda}g,g}\rangle_{\pi}\nonumber\\
&\ge&\left(\int_{\sigma(Q)}\lambda^{k}\langle{d\,E_{\lambda}g,g}\rangle_{\pi}\right)^{\frac{n}{k}}=\langle{Q^{k}g,g}\rangle^{\frac{n}{k}}_{\pi}.
\end{eqnarray}
This together with (\ref{eins}) and (\ref{laterneed}) yields
\begin{eqnarray}\label{drei}
\langle{Q^{k}g,g}\rangle_{\pi}&\le&\langle{Q^{n}g,g}\rangle^\frac{k}{n}_{\pi}\le\left(\langle{P^{\ast^n}P^n g,g}\rangle_{\pi}+|\langle{(P^{\ast^n}P^n -Q^n)g,g}\rangle_{\pi}|\right)^\frac{k}{n}\nonumber\\
&\le& \max\left(2\langle{P^{\ast^n}P^n g,g}\rangle_{\pi},2|\langle{(P^{\ast^n}P^n -Q^n)g,g}\rangle_{\pi}| \right)^{\frac{k}{n}}\nonumber\\
&\le& \max\left(\left(2l_{g}(1-k_{P^{\ast^n}P^n})\right)^{\frac{k}{n}},\left(2||(P^{\ast^n}P^n -Q^n)g||_{2}\right)^{\frac{k}{n}}\right)
\end{eqnarray}
Now with $n\rightarrow\infty$ we obtain for all $k\in\mathbb{N}$
\begin{equation}\label{hutti}
\langle{Q^{k}g,g}\rangle_{\pi}\le \max\left(\lim_{n\rightarrow\infty}\left(1-k_{P^{\ast^n}P^n}\right)^{\frac{k}{n}}, q^k\right).
\end{equation}
At this point note the small but crucial difference between (\ref{hutti}) and (\ref{laterneed}). By taking advantage of the reversible structure of $Q$, we get rid of the therm $l_g$. Now 
the result follows since $f\in L^{2}_{0,1}(\pi)$ can be approximate arbitrary well by indicator functions $g\in L^{2}_{0,1}(\pi)$.
\finishproof

As an immediate consequence we obtain

\begin{cor}\label{isorhocor}
Assume that the operator $P$ associated to the Markov chain $\xi_1,\xi_2,\ldots$ is normal. Then the spectral radius $\rho$
is given by
\begin{equation}\label{dernormalefall}
\rho=\lim_{n\rightarrow\infty} \left(\sqrt{1-k_{P^{\ast^n}P^n}}\right)^{\frac{1}{n}}.
\end{equation}
Especially, for reversible Markov chains this yields
\begin{equation}
\rho=\lim_{n\rightarrow\infty} \left(\sqrt{1-k_{2n}}\right)^{\frac{1}{n}}.
\end{equation}
Moreover, if the Markov chain is also positive, we obtain
\begin{equation}
\rho=\lim_{n\rightarrow\infty} \left(1-k_{n}\right)^{\frac{1}{n}}.
\end{equation}
\end{cor}
\proof
We see that for $Q:=P^{\ast}P$, (\ref{q}) is satisfied with $q=0$. Now the claim follows from Theorem \ref{assumdiscuss}
and (\ref{evenequality2}). The second assertion is trivial, since $P^{\ast}P=P²$ in the self adjoint case. The last 
assertion follows from the fact that can do the proof of Theorem \ref{assumdiscuss} directly using $P$ instead of $Q$.
\finishproof
Let us interpret this results in terms of spectral measures: Denote $\mu_{f}(d\lambda)=d\langle{E_{\lambda}f,f}\rangle$ the spectral measure according to the Markov operator $P$ and $f\in L^{2}_{0,1}(\pi)$. 
Then the Rayleigh-Ritz principle for reversible Markov chains implies that there exists a sequence of functions
$(f_i)_{i\in\mathbb{N}}\subset L^{2}_{0,1}(\pi)$ such that for all $\epsilon>0$ we have
\begin{equation}
\lim_{i\rightarrow\infty}\mu_{f_i}(B_{\rho}(\epsilon))=1 ,
\end{equation}
where $B_{\rho}(\epsilon):=\{y\in\mathbb{R}:|y-\rho|\le\epsilon\}$.
Let us consider a certain, non-linear subspace of $L^{2}_{0,1}(\pi)$, namely the step-functions of the following kind:
\begin{equation}\label{Maul}
f_{A}=\sqrt{\pi(A)\pi(A^{c})}\left(\frac{1}{\pi(A)}1_{A}-\frac{1}{\pi(A^{c})}1_{A^{c}}\right), \,\,A\in\mathcal{F}.
\end{equation}
It can be easily checked that 
\begin{equation}\label{Maul2}
E_{\pi}f_{A}=0,\,\,\,||f_{A}||_2 =1.
\end{equation}
Moreover, we have the following identity:
\begin{equation}
1-k_{n}(A)=\langle{f_{A},P^{n}f_{A}}\rangle_{\pi}=\int_{-1}^{1}\lambda^n\mu_{f_{A}}(d\lambda). 
\end{equation}
Since $\rho=\lim_{n\rightarrow\infty}(1-k_{2n})^{\frac{1}{2n}}$, it follows that there exists a sequence 
$(A_i)_{i\in\mathbb{N}}\subset\mathcal{F}$ such that for all $\epsilon>0$ we have that 
\begin{equation}
\lim_{i\rightarrow\infty}\mu_{f_{A_i}}(B_{\rho}(\epsilon))>0 ,
\end{equation}
In comparison to the general Rayleigh-Ritz-principle, which says that the family of spectral measures $\mu_{f},f\in
L^{2}_{0,1}(\pi)$ contains measures that are in a certain sense arbitrary close to the dirac-measure at $\rho$, we see that the
family of spectral measures induced by the family of simple step functions defined by (\ref{Maul}) contains for every $\epsilon>0$ at least some measures that put a positive amount of mass to the ball $B_{\rho}(\epsilon)$.

\section{Applications}\label{apply}
In this section we will present some applications of the developed theory. The first applications should indicate the theoretical value of Theorem \ref{assumdiscuss}. Later we will see how to use Theorem \ref{assumdiscuss} in order to obtain precise bounds for the $M/M/1$-queue, the hypercube
and the star.
Let us start with a partially generalization of a Theorem due to Roberts and Tweedie \cite{roberts1}:
\begin{theorem}
Assume that $\xi_1,\xi_2,\ldots$ is a geometrically ergodic and normal Markov chain with (ORGC) $\delta$. Then
$\xi_1,\xi_2,\ldots$ has the (SGP) and the spectral radius $\rho$ can be estimated by 
\[\rho\le\delta\]
\end{theorem}
\proof

From (\ref{zwei}) it follows that  
\begin{equation}\label{doeblin4}
(1-k_{P^{\ast^l}P^l}(A))^{\frac{1}{2l}}\le(1-k_{P^{\ast^n}P^n}(A))^{\frac{1}{2n}}. 
\end{equation}
In (\ref{adjungrate}) we will show that 
\begin{equation}\label{doeblin3}
1-k_{P^{\ast^n}P^n}(A)=\frac{1}{\pi(A)\pi(A^{c})}\int_{\Omega}(p^{n}(x,A^{c})-\pi(A^{c}))^{2}\pi(dx).
\end{equation}
Now, using (\ref{doeblin4}), (\ref{doeblin3}) and Theorem \ref{habbatu} it follows that
\begin{eqnarray}
(1-k_{P^{\ast^l}P^l}(A))^{\frac{1}{2l}}&\le&\left(\frac{1}{\pi(A)\pi(A^{c})}\int_{\Omega}(p^{n}(x,A^{c})-\pi(A^{c}))^{2}\pi(dx)\right)^{\frac{1}{2n}}\nonumber\\.
&\stackrel{}{\le}& \left(\frac{2}{\pi(A)\pi(A^{c})}\right)^{\frac{1}{2n}}\delta 
\end{eqnarray}
Now let first $n\rightarrow\infty$, then take the infimum over all $A\in\mathcal{F}$ and finally let $l\rightarrow\infty$. Now the claim follows from Theorem \ref{assumdiscuss}.
\finishproof

Note that this result is known for reversible Markov chains \cite{roberts1}. The existing proof for that case is more difficult and requires results that have been appeared in \cite{roberts2}.  \\

Let us improve the results obtained in Corollary \ref{hidoe}.

\begin{cor}\label{doeblin_revisted}
Assume that the Markov chain $\xi_1,\xi_2,\ldots$ satisfies (\ref{geometrischergodisch}) with \\
$C(x)=C$ independent of $x$ and that additionally $P^{\ast}$ has a representation of the form
\begin{equation}
P^{\ast}f(x):=\int_{\Omega}f(y)p^{\ast}(x,dy),
\end{equation}
such that $p^{\ast}(x,\cdot)$ is a transition-probability kernel.  
Then $P$ has an $L^{2}(\pi)$-spectral gap and the spectral 
radius $\rho$ can be estimated by $\sqrt{\delta}$.
\end{cor}
\proof

Putting (\ref{doeblin3}) and (\ref{doeblin4}) together we obtain for $l_{1}<l_{2}\in\mathbb{N}$ that
\begin{eqnarray}
(1-k_{(P^{\ast^{n}}P^{n})^{2l_1}}(A))^{\frac{1}{2l_1}}&\le& (1-k_{(P^{\ast^{n}}P^{n})^{2l_2}}(A))^{\frac{1}{2l_2}}\nonumber\\
&\le&\left(\frac{1}{\pi(A)\pi(A^{c}}\right)^{\frac{1}{2l_2}}\nonumber\\
&&\,\,\,\times\left(\int_{\Omega}\left((P^{\ast^{n}}P^{n})^{2l_2}(1_{A}-\pi(A))\right)^2(x)\pi(dx)\right)^{\frac{1}{2l_2}}\nonumber\\
&\le&\left(\frac{1}{\pi(A)\pi(A^{c}}\right)^{\frac{1}{2l_2}}\nonumber\\
&&\left(\sup_{x\in\Omega}||p^{\ast^{n}}(x,\cdot)-\pi||_{TV}^{2l_2} \sup_{x\in\Omega}||p^{n}(x,\cdot)-\pi||_{TV}^{2l_2}\right)^{\frac{1}{2l_2}}\nonumber\\
\end{eqnarray}
Now let $l_{2}\rightarrow\infty$ and take the supremum over all $A\in\mathcal{F}$ we obtain
\begin{equation}
(1-k_{(P^{\ast^{n}}P^{n})^{2l_1}})^{\frac{1}{2l_1}}\le 2C\delta^{n}.
\end{equation}
Since $P^{\ast^{n}}P^{n}$ is selfadjoint and positive, we can apply Corollary \ref{isorhocor}. For $l_1\rightarrow\infty$ we get
\begin{equation}
\rho(P^{\ast^{n}}P^{n})\le 2C\delta^{n},
\end{equation}
where $\rho(P^{\ast^{n}}P^{n})$ denotes the spectral radius of $P^{\ast^{n}}P^{n}$ on $L^{2}_{0}$. Since
\begin{equation}
\rho(P^{\ast^{n}}P^{n})=\sup_{f\in L^{2}_0}\langle{P^{\ast^{n}}P^{n}f,f}\rangle_{\pi}=||P^n||^{2}_{L^{2}_0},
\end{equation}
with $||P^n||_{L^{2}_{0}}$ the operator-norm, we finally obtain, by taking the $2n-th$-square root and let $n\rightarrow\infty$, that
\begin{equation}
\rho\le \sqrt{\delta}.
\end{equation}
\finishproof

This result is well-known and can be alternatively obtained by applying Riesz-Thorin's-interpolation theorem (see e.g. \cite{werner}) in
the following way. By Riesz-Thorins-theorem it follows that 
\begin{equation}\label{thorin1}
||f||_{L^{2}_0}\le||f||_{L^{\infty}_0}^{\frac{1}{2}}||f||_{L^{1}_0}^{\frac{1}{2}}.
\end{equation}
In \cite{wuebkerdiss}, Proposition 1, the following inequality was shown:
\begin{equation}\label{thorin2}
\sup_{x\in\Omega}||p^{n+1}(x,\cdot)-\pi||_{TV}\le ||P^{n}-P_1||_{\infty}\le\sup_{x\in\Omega}||p^{n}(x,\cdot)-\pi||_{TV}.
\end{equation}
Here, $||\cdot||_{\infty}$ denotes the supremum-norm over all bounded, measurable complex valued functions. 
Now from (\ref{thorin1})
and (\ref{thorin2}) also we also obtain the result of Corollary \ref{doeblin_revisted}. But note that the proof of
Corollary \ref{doeblin_revisted} actually yields more. In the proof we used the following trivial estimate:
\[\sup_{x\in\Omega}||p^{\ast^{n}}(x,\cdot)-\pi||_{TV}\le 2.\]
But if we have that $p^{\ast^{n}}(x,\cdot)$ is also uniformly ergodic, which is for finite Markov chains always
the case, we immideately obtain
\begin{cor}
Let $\xi_1,\xi_2,\ldots$ be a Markov chain as in Corollary \ref{doeblin_revisted} and assume, that additionally,
there exists a constant $C^{\ast}>0$ and $\delta^{\ast^n}<1$ such that
\begin{equation}
\sup_{x\in\Omega}||p^{\ast^{n}}(x,\cdot)-\pi ||_{TV}\le C^{\ast}\delta^{\ast^n}.
\end{equation}
Then
\begin{equation}
\rho\le\sqrt{\delta\delta^{\ast}}.
\end{equation}
\end{cor}
\proof
See the proof of Corollary 6
\finishproof

Now the theorems will be applied to obtain explicit bounds for some 
Markov chains.

\begin{example}[M/M/1-Queue]
Let us consider the $M/M/1$ queuing system with state space $\mathbb{N}$ and reflecting boundary at one. More precisely, the non
zero elements of $P=(p_{i,j})_{(i,j)\in\mathbb{N}x\mathbb{N}}$ are
\begin{eqnarray}
\,\,\,\quad\quad p_{i,i+1}&=&q,\,\,\,\forall i\in\mathbb{N}\nonumber\\
\,\,\,\quad\quad p_{i,i-1}&=&p,\,\,\, \mathbb{N}\ni i\ge 2\nonumber\\
\,\,\,\quad\quad p_{1,1}&=&p,
\end{eqnarray}
where $p>\frac{1}{2}$ in order to guarantee that the associated Markov chain is positive recurrent. This example was studied in the book 
of Feller (\cite{feller}, page 436-438). From that one can deduce (for example by using that spectral radius=(ORGC) for reversible Markov chains) that the $L^{2}(\pi)$-spectral radius $\rho$ of this chain is given by 
\begin{equation}\label{optimal_bound}
\rho=2\sqrt{pq}.
\end{equation}  
Now we will proceed by showing how to derive this bound by using isoperimetric constants. For $s\in\mathbb{N}$ let
\[A_{s}:=\{s,s+1,s+2,\ldots\}\]
The invariant measure $\pi$ can be easily calculated to be
\[\pi(i)=\frac{2p-1}{p}\left(\frac{q}{p}\right)^{i-1},\,\,\,\forall i\in\mathbb{N}.\]
It follows that for all $s\in\mathbb{N}$ we have 
\[\pi(A_{s})=\left(\frac{q}{p}\right)^{s-1}.\]
This yields
\begin{eqnarray}\label{mm1}
k_{2n}(A_s)&=&\frac{1}{\left(\frac{q}{p}\right)^{s-1}\left(1-\left(\frac{q}{p}\right)^{s-1}\right)}\sum_{i=s}^{\infty}p^{2n}(i,\{1,2,\ldots,s-1\})\frac{p-q}{p}\left(\frac{q}{p}\right)^{i-1}\nonumber\\
&=&\frac{1}{\left(1-\left(\frac{q}{p}\right)^{s-1}\right)}\frac{p-q}{p}\sum_{i=0}^{\infty}p^{2n}(s+i,\{1,2,\ldots,s-1\})\left(\frac{q}{p}\right)^{i}.
\end{eqnarray}
The interesting observation now is that for $s>n$, (\ref{mm1}) is decreasing in $s$, and one might hope that $\lim_{s\rightarrow\infty}k_{n}(A_s)=k_n$.
We have
\begin{equation}\label{mm12}
\lim_{s\rightarrow\infty}k_{2n}(A_s)=\frac{p-q}{p}\sum_{i=0}^{\infty}p^{2n}(s+i,\{1,2,\ldots,s-1\})\left(\frac{q}{p}\right)^{i}.
\end{equation}
Note that $p^{2n}(s+i,\{1,2,\ldots,s-1\})$ does not depend on $s$ for $s>n$. Now (\ref{mm12}) implies that
\begin{equation}\label{mm13}
1-\lim_{s\rightarrow\infty}k_{2n}(A_s)=\frac{p-q}{p}\sum_{i=0}^{\infty}p^{2n}(s+i,\{s,s+1,s+2,\ldots\})\left(\frac{q}{p}\right)^{i}.
\end{equation}
Since for $s>n$ we have
\begin{equation}\label{mm14}
p^{2n}(s+i,\{s,s+1,s+2,\ldots\})=\sum_{j=0}^{i}{2n\choose[\frac{2n+j}{2}]}p^{[\frac{2n+j}{2}]}q^{2n-[\frac{2n+j}{2}]}.
\end{equation}
Now we apply (\ref{mm13}) to (\ref{mm14}) and obtain
\begin{eqnarray}
1-\lim_{s\rightarrow\infty}k_{2n}(A_s)&=&\frac{p-q}{p}\sum_{i=0}^{2n}\sum_{j=0}^{i}{2n\choose[\frac{2n+j}{2}]}p^{[\frac{2n+j}{2}]}q^{2n-[\frac{2n+j}{2}]}\left(\frac{q}{p}\right)^{i}\label{notrough}\\
&\ge&\frac{p-q}{p}{2n\choose n}p^{n}q^{n}.\label{veryrough}
\end{eqnarray}
Using Stirling's approximation and apply Corollary \ref{isorhocor} one obtains
\begin{equation}\label{mm15}
\rho=\lim_{n\rightarrow\infty}(1-k_{2n})^{\frac{1}{2n}}\ge \lim_{n\rightarrow\infty}(1-\lim_{s\rightarrow\infty}k_{2n}(A_s))^{\frac{1}{2n}}\ge 2\sqrt{pq}.
\end{equation}
On the other hand, the right hand side in (\ref{notrough}) can be estimated from above by
\begin{equation}
\frac{p-q}{p}{2n\choose n}p^n q^n\sum_{i=0}^{2n}\sum_{j=0}^{i}\left(\frac{q}{p}\right)^{i-[\frac{j}{2}]}<\frac{p-q}{p}4n^2 {2n\choose n}p^n q^n.
\end{equation}
Again, applying Stirling's formula yields 
\begin{equation}
\lim_{n\rightarrow\infty}(1-\lim_{s\rightarrow\infty}k_{2n}(A_s))^{\frac{1}{2n}}\le 2\sqrt{pq}.
\end{equation}
This together with (\ref{mm15}) gives
\begin{equation}\label{mm16}
(1-\lim_{s\rightarrow\infty}k_{2n}(A_s))^{\frac{1}{2n}}= 2\sqrt{pq}.
\end{equation}
We have just shown that $\rho\ge 2\sqrt{pq}$ and if we also have that
\begin{equation}\label{wanne}
\left(1-\lim_{s\rightarrow\infty}k_{2n}(A_s)\right)^{\frac{1}{2n}}=(1-k_{2n})^{\frac{1}{2n}},
\end{equation}
then $\rho=2\sqrt{pq}$. On the opposite, since we know already that $\rho=2\sqrt{pq}$, it follows that (\ref{wanne}) holds
\end{example}
Note that the estimate in (\ref{veryrough}) is very rough but still suffices to obtain a sharp lower bound for the spectral radius $\rho$.\\
Let us recapitulate the procedure. In order to obtain a "good" lower bound for the spectral radius $\rho$, we choose from our intuition (or our
knowledge) a sequence $(A_s)_{s\in\mathbb{N}}$ such that it is supposed that 
\begin{equation}\label{mm17}
\lim_{s\rightarrow\infty}k_{n}(A_s)\approx k_n.
\end{equation} 
Next we calculate, if possible, the left hand side in (\ref{mm17}) and obtain a lower bound. If in fact (\ref{mm17}) holds
true in a certain sense, the archived lower bound for the spectral radius is sharp. Unfortunately, in most cases it is rather difficult to 
verify (\ref{mm17}), so one depends on a "`good"' feeling for the right sequence $(A_s)_{s\in\mathbb{N}}$. Nevertheless, we will give two 
more examples where the intuitive choice of $(A_s)_{s\in\mathbb{N}}$ yield sharp lower bounds for the spectral radius $\rho$.

\begin{example}[Lazy random walk on hypercube]
We will show how we can apply Theorem \ref{assumdiscuss} to obtain a sharp upper estimate of the spectral gap of the lazy random walk on the hypercube by
analyzing the spectral properties of a much simpler Markov chain.
Let us consider the random walk on the hypercube with the set $V$ of vertices given by 
\[V:=\{x=(x_1,x_2,\ldots,x_n): x_i\in\{0,1\}, i\in\{1,2,\ldots,n\}\}.\]
The transition probabilities $p(\cdot,\cdot)$ of the random walk are given by
\[p(x,y)=\left\{\begin{array}{r@{\quad:\quad}l}
\frac{1}{2}& x=y\\
\frac{1}{2n}& ||x-y||_2=1\\
0& \mbox{else}\\
\end{array}\right..\]
Define 
\[A=\{x\in V: \,x_1=0\}.\]
For all $x\in A$, $m\in\Nset$ we have 
\begin{equation}\label{seclareigen}
p^{m}(x,A^{c})=(1,0)
Q^{m}{0 \choose 1},
\end{equation} 
where $Q:=\left(\begin{array}{c@{\quad}c}
1-\frac{1}{2n} & \frac{1}{2n}\\
\frac{1}{2n} & 	1-\frac{1}{2n}
\end{array}
\right)$.
An easy calculation shows 
\[p^{m}(x,A^{c})=\frac{1}{2}-\frac{1}{2}(1-\frac{1}{n})^{m}.\]

From this we obtain
\[\lim_{m\rightarrow\infty}(1-k_{2m}(A))^{\frac{1}{2m}}=\lim_{m\rightarrow\infty}\left(\left(1-\frac{1}{n}\right)^{2m}\right)^{\frac{1}{2m}}=1-\frac{1}{n},\]
and therefore
\[\rho\ge 1-\frac{1}{n}.\]
In fact, it is known that $\rho=1-\frac{1}{n}$, so the bound is sharp.  
\end{example}

\begin{example}[The star]
Consider the following Markov chain with state space $\{1,2,\ldots,n\}$ Let $a_i \in (0,1),i \in \{1,2,\ldots,n\}$ such that
\[\sum_{i=1}^{n}a_i=1.\]
Let the transition probabilities be given by
\begin{equation}
P=\left(
\begin{array}{c c c c c}
a_1&a_2&a_3&\ldots &a_n\\
1&0&0&\ldots&0\\
1&0&0&\ldots&0\\
\vdots&\vdots&\vdots&\ldots&0\\
1&0&0&\ldots&0
\end{array}
\right).
\end{equation}
The invariant distribution $\pi$ of this Markov chain can be easily verified to be 
\[\pi(1)=\frac{1}{2-a_1},\,\,\,\pi(i)=\frac{a_i}{2-a_1},\,\,i\in\{2,3,\ldots,n\}.\]
It can be checked that the chain is reversible, so Theorem \ref{assumdiscuss} can be 
applied. Since everything flows out of the set $A:=\{2,3,\ldots,n\}$ in one step, one expects that \[\lim_{n\rightarrow\infty}(1-k_{2n}(A))^{1/2n}=\lim_{n\rightarrow\infty}(1-k_{2n})^{1/2n}.\] 
Again, the left hand side
can be calculated as in the example before with the Matrix 
\[Q:=\left(\begin{array}{c@{\quad}c}
0 & 1\\
a_1 & 	1-a_1
\end{array}
\right)\]
The second largest eigenvalue of this Matrix equals
$\,-a_1$. This yields
\[\rho\ge a_1.\]
Actually, it can be shown that $\rho=a_1$, so the bound is sharp. Again, we have just seen that 
Theorem \ref{assumdiscuss} can be used to reduce the complexity of the original chain by analyzing a much simpler chain. 
The following observation may be interesting: The special structure of the star yields to the fact that independent of the 
chosen set $A$ it can be shown that $\lim_{n\rightarrow\infty}(1-k_{2n}(A))^{1/2n}=\lim_{n\rightarrow\infty}(1-k_{2n})^{1/2n}$.
This is in some sense good and bad. It is good because it shows that even an arbitrary chosen set $A$ may yield a sharp lower bound
for the spectral radius of $P-P_1$. It is bad because it shows that it might be impossible to recover all or most eigenvalues of the 
Markov chain by computing $\lim_{n\rightarrow\infty}(1-k_{2n}(A))^{1/2n} $ for different sets $A$. In this example, the eigenvalue $0$ with
multiplicity $n-2$ is not recognized.

\end{example}

\section{Isoperimetric constants and mixing}\label{mixing1}
This section is included to interpret the obtained results in terms of mixing.
Let $(X_n)_{n\in\mathbb{N}}$ denote a family of random variables on a probability space $(\Omega,\mathcal{F},P)$ and let 
$\mathcal{F}_{m}:=\sigma(X_1,X_2,\ldots,X_m)$ and $\mathcal{F}^{n}:=\sigma(X_{n},X_{n+1},\ldots)$. The sequence $(X_n)_{n\in\mathbb{N}}$
is called strongly mixing ($\alpha$-mixing) if $\alpha(n)\rightarrow 0$, where 
\[\alpha(n):=\sup_{k\ge 1}\sup_{A\in \mathcal{F}_{k},B\in\mathcal{F}^{k+n}}|P(A\cap B)-P(A)P(B)|,\]
uniform mixing ($\phi$-mixing) if $\phi(n)\rightarrow 0$ with
\begin{equation}\label{phimixing}
\phi(n):=\sup_{k\ge 1}\sup_{A\in \mathcal{F}_{k},B\in\mathcal{F}^{k+n}}\frac{|P(A\cap B)-P(A)P(B)|}{P(A)},
\end{equation}
and asymptotically uncorrelated ($\rho$-mixing) if $\rho(n)\rightarrow 0$ with
\[\rho(n):=\sup\{corr(U,V), U\in L^{2}(\mathcal{F}_{k}), V\in L^{2}(\mathcal{F}^{k+n}).\]
It is well-known that $\phi$-mixing implies $\rho$-mixing and $\rho$-mixing implies $\alpha$-mixing (see e.g. \cite{jones}). In fact,
$\rho$-mixing is equivalent to the (SGP) and $\phi$-mixing is equivalent to uniform ergodicity (\cite{jones}). The question we will answer is how to compare $\rho$-mixing with $\phi$-mixing from their definitions. We want to have a mixing definition which is equivalent to $\rho$-mixing, but comparable with
$\phi$-mixing in a sense as geometric ergodicity is comparable to uniform ergodicity.\\  
During the last two sections we have seen that the (SGP) is equivalent to
\begin{equation}\label{mixing}
1-k_n\le C\delta^n,
\end{equation}
for some $C>0$ and $\delta<1$. But this is equivalent to
\begin{equation}\label{mixingtwo}
\sup_{k\ge 1}\sup_{A_{k}\in \sigma(\xi_k),A_{n+k}^{c}\in\sigma(\xi_{n+k})}\frac{|P(A_k\cap A_{k+n}^{c} )-P(A_{k})P(A_{k+n}^{c})|}{P(A_{k})}\le\tilde{C}\delta^{n},
\end{equation}
for some $\tilde{C}>0$, $\delta<1$ as in (\ref{mixing}) and 
\[A_k:=\Omega_1\times\Omega_2\times\ldots\Omega_{k-1}\times A\times\Omega_{k+1}\times\ldots\,\, \mbox{with }\Omega_i =\Omega.\]
So (\ref{mixingtwo}) can be seen as a suitable way to weaken 
(\ref{phimixing}) in such a way that an equivalent condition to the (SGP) is obtained. Moreover, in the reversible case
the optimal chosen $\delta$ on the right hand side in (\ref{mixingtwo}) equals the spectral radius of $P-P_1$ on $L^{2}(\pi)$.

\section{Proofs}\label{proofs}
\proof
In (\ref{Maul}) we considered the function
\begin{equation}\label{fnot}
f_{0}=\sqrt{\pi(A)\pi(A^{c})}\left(\frac{1}{\pi(A)}1_{A}-\frac{1}{\pi(A^{c})}1_{A^{c}}\right), \,\,A\in\mathcal{F}.
\end{equation}
We have already seen in (\ref{Maul2}) that
\[E_{\pi}f_{0}=0,\,\,\,||f_{0}||_2 =1.\]
From the definition of $f_{0}$ we calculate
\begin{eqnarray}
||Pf_{0}||_{2}^{2}&=&\pi(A)\pi(A^{c})\int_{\Omega}\left(\frac{p(x,A)}{\pi(A)}-\frac{p(x,A^{c})}{\pi(A^{c})}\right)^{2}\pi(dx)\nonumber\\
&=&\pi(A)\pi(A^{c})\left(\frac{1}{\pi(A)^2}\int_{\Omega}p(x,A)^2\pi(dx)+\frac{1}{\pi(A^{c})^2}\int_{\Omega}p(x,A^{c})^2\pi(dx)\right.\nonumber\\
&& \left.-2\frac{1}{\pi(A)\pi(A^{c})}\int_{\Omega}p(x,A)p(x,A^{c})\pi(dx)\right).\label{einsetzen}
\end{eqnarray}
We have
\begin{eqnarray}
\frac{1}{\pi(A)^2}\int_{\Omega}p(x,A)^2\pi(dx)&=&\frac{1}{\pi(A)^2}\left(\int_{A}p(x,A)^2\pi(dx)+\int_{A^{c}}p(x,A)^2\pi(dx)\right)\nonumber\\
&=& \frac{1}{\pi(A)^2} \left(\int_{A}(1-p(x,A^{c}))^2\pi(dx)+\int_{A^{c}}p(x,A)^2\pi(dx)\right)\nonumber\\
&=&\frac{1}{\pi(A)^2}\left(\pi(A)-2\pi(A)\pi(A^{c})k(A)+\int_{A}p(x,A^{c})^{2}\pi(dx)\right.\nonumber\\
&&+\left.\int_{A^{c}}p(x,A)^{2}\pi(dx)\right)\nonumber\\
&\stackrel{Jensen}{\ge}&\frac{1}{\pi(A)}(1-\pi(A^{c})k(A)(2-k(A)))\nonumber.
\end{eqnarray}
Using $k(A)=k(A^{c})$ one obtains
\begin{equation}
\frac{1}{\pi(A^{c})^2}\int_{\Omega}p(x,A^{c})^2\pi(dx)\ge \frac{1}{\pi(A^{c})}(1-\pi(A)k(A)(2-k(A))).
\end{equation}
Inserting the last three inequalities in (\ref{einsetzen}) yields
\begin{eqnarray}
||Pf_{0}||_{2}^{2}&\ge&\pi(A^{c})-\pi(A^{c})^{2}k(A)(2-k(A))+\pi(A)-\pi(A)^{2}k(A)(2-k(A))\nonumber\\
&&-2\pi(A)\pi(A^{c})k_{P^{\ast}P}(A)\nonumber\\
&=&(1-k(A))^{2}+2\pi(A)\pi(A^{c})\left((1-k_{P^{\ast}P}(A))-(1-k(A))^{2}\right)\nonumber\\
&\ge&\frac{1}{2}(1-k(A))^{2}+2\pi(A)\pi(A^{c})(1-k_{P^{\ast}P}(A)).
\end{eqnarray}
The same calculations can be done with $P^{n}$ instead of $P$. So the claim follows.
\finishproof


Now let us proof Theorem \ref{spaet}

\proof
Since (\ref{wichtigevoraussetzung})-which was established to be necessary for the existence of a spectral gap (Proposition \ref{noetig}-implies (\ref{zwo}), we already 
know that (\ref{zwo}) is necessary. It remains to show that this condition is also sufficient.\\
So let us assume that (\ref{zwo}) holds. Then there exists $\epsilon>0$ such that

\begin{equation}\label{aerger1}
\Re(\sigma(P^{n}))\le 1-\frac{\kappa}{8}\epsilon^2\,\,\,\forall n\le\left[\frac{2\pi}{\arccos(1-\frac{\kappa}{16}\epsilon^{2})}\right]+1.
\end{equation}
On the other hand, from the spectral mapping theorem we know that
\begin{equation}\label{aerger2}
\sigma(P^{n})=(\sigma(P))^{n}.
\end{equation}
Now let $z\in\sigma(P)$. Using polar coordinates, $z$ has the representation
\[z=z(r,\phi)=re^{i\phi}=r\cos(\phi)+ir\sin(\phi),\,\,r\in[0,1],\,\,\phi\in[-\pi,\pi).\]
From (\ref{aerger1}) and (\ref{aerger2}) we obtain
\[r^{n}\cos n\phi\le 1-\frac{\kappa}{8}\epsilon^{2},\,\,\forall n\le\left[\frac{2\pi}{\arccos(1-\frac{\kappa}{16}\epsilon^{2})}\right]+1.\]
Now choose $\phi_{1}=\arccos(1-\frac{\kappa}{16}\epsilon^{2})\cap (0,\pi)$. For all $\phi\in(-\phi_1,\phi_1)$ and $z(r,\phi)\in\sigma(P)$, we obtain: $r\le 1-\frac{\kappa}{16}\epsilon^{2}$.
For those $z(r,\phi)$ with $\phi\not\in(-\phi_1,\phi_1)$ there exists $n=n(\phi)\le n_{0}(\phi_1):=\left[\frac{2\pi}{\phi_1}\right]+1$, such that $n\phi[mod 2\pi]\in(-\phi_1,\phi_1)$.
This together with (\ref{aerger2}) yields: \[r\le(1-\frac{\kappa}{16}\epsilon^{2})^{1/([\frac{2\pi}{\arccos(1-\frac{\kappa}{16}\epsilon^{2})}]+1)}.\]
\finishproof

We should now prove Lemma \ref{dasbraucheichnoch}.\\
\proof
For $i\in\{1,2,\ldots,n_0\}$ consider
\begin{eqnarray}
\pi(A)\pi(A^{c})k_{n}(A)&=&\int_{A}p^{n}(x,A^{c})\pi(dx)=\int_{A}\int_{A}p(x,dy)p^{n-1}(y,A^{c})\pi(dx)\nonumber\\
&+&\int_{A}\pi(dx)\int_{A^{c}}p(x,dy)p^{n-1}(y,A^{c})\nonumber\\
&\le&(k_{n-1}(A)+k(A))\pi(A)\pi(A^{c})\le\ldots\le n\,k_{1}(A)\pi(A)\pi(A^{c})\nonumber.
\end{eqnarray}

For the other
$k_i$, $i\in\{2,\ldots,n_{0}-1\}$ we have
\[\int_{A}p^{i}(x,A^{c})\ge \frac{1}{2}\int_{A}p^{2i}(x,A^{c})\ge\ldots \ge\frac{1}{2^{[\log_{2} n_0]}}\int_{A}p^{i2^{[\log_{2} n_0]}}(x,A^{c}).\]
Since
\[i2^{[\log_{2} n_0]}\ge n_{0}\,\,\forall i\in\{2,\ldots,n_{0}-1\}\] 
and
\[2^{[\log_{2} n_0]}\le n_{0},\] 
the claim follows.
\finishproof

Let us present the proof of Lemma \ref{monotonie}:

\proof
First let us proof the second statement:
\begin{eqnarray}
k_{P^{\ast^{n}}P^{n}}(A)&=&\frac{1}{\pi(A)\pi(A^{c})}\int_{A}P^{\ast^{n}}P^{n}1_{A^{c}}\pi(dx)\nonumber\\
&=&\frac{1}{\pi(A)\pi(A^{c})}\int_{\Omega}p^{n}(x,A)p^{n}(x,A^{c})\pi(dx)\nonumber\\
&=&\frac{1}{\pi(A)}-\frac{1}{\pi(A)\pi(A^{c})}\int_{\Omega}(p^{n}(x,A^{c}))^{2}\pi(dx)\nonumber\\
&\stackrel{Jensen}{\le}&\frac{1}{\pi(A)}-\frac{\pi(A^{c})}{\pi(A)}=1.\nonumber
\end{eqnarray}
The first assertion follows from
\begin{eqnarray}
k_{P^{\ast^{n+1}}P^{n+1}}(A)&=&\frac{1}{\pi(A)}-\frac{1}{\pi(A)\pi(A^{c})}\int_{\Omega}(p^{n+1}(x,A^{c}))^{2}\pi(dx)\nonumber\\
&\stackrel{P\,\, Markov}{\ge}&\frac{1}{\pi(A)}-\frac{1}{\pi(A)\pi(A^{c})}\int_{\Omega}P(p^{n}(x,A^{c}))^{2}\pi(dx)=k_{P^{\ast^{n}}P^{n}}(A)\nonumber.
\end{eqnarray}
\finishproof

We should now present the proof of Lemma \ref{bani}

\proof
Without loss of generality we may assume that $\pi(A)\le\frac{1}{2}$. On the one hand we have
\begin{eqnarray}
1-k_{n}(A)&=&1-\frac{1}{\pi(A)\pi(A^{c})}\int_{A}p^{n}(x,A^{c})\pi(dx)\nonumber\\
&=&\int_{A}(\pi(A^{c})-p^{n}(x,A^{c}))\frac{\pi(dx)}{\pi(A)\pi(A^{c})}\nonumber\\
&\stackrel{C.S.}{\le}&\left(\int_{A}(\pi(A^{c})-p^{n}(x,A^{c}))^{2}\frac{\pi(dx)}{\pi(A)\pi(A^{c})}\right)^{\frac{1}{2}}\frac{1}{\sqrt{\pi(A^{c})}}\nonumber\\
&\le&\sqrt{2}\sqrt{\frac{1}{\pi(A)\pi(A^{c})}\int_{\Omega}(\pi(A^{c})-p^{n}(x,A^{c}))^{2}\pi(dx)}.\label{vergl.0}
\end{eqnarray}
On the other hand we get
\begin{eqnarray}\label{adjungrate}
1-k_{P^{\ast^n}P^n}(A)&=&1-\frac{1}{\pi(A)\pi(A^{c})}\int_{\Omega}p^{n}(x,A)p^{n}(x,A^{c})\pi(dx)\nonumber\\
&=&1-\frac{1}{\pi(A)}+\frac{1}{\pi(A)\pi(A^{c})}\int_{\Omega}p^{n}(x,A^{c})^{2}\pi(dx)\nonumber\\
&=&\frac{1}{\pi(A)\pi(A^{c})}\int_{\Omega}(p^{n}(x,A^{c})^2-\pi(A^{c})^{2})\pi(dx)\nonumber\\
&=&\frac{1}{\pi(A)\pi(A^{c})}\int_{\Omega}(p^{n}(x,A^{c})-\pi(A^{c}))^{2}\pi(dx).\label{vergl.}
\end{eqnarray}
Reinsert (\ref{vergl.}) in (\ref{vergl.0}), we achieve
\begin{equation}
1-k_{n}(A)\le\sqrt{2}\sqrt{1-k_{P^{\ast^n}P^n}(A)}.
\end{equation}
Now taking the supremum first on the right and thereafter on the left hand side, the claim follows.
\finishproof

\section{Acknowledgments}
The author thanks Wolfgang Stadje and Zakhar Kabluchko for reading the paper.


\end{document}